\documentclass{article}                            
\usepackage{amsmath,amssymb,citesort,enumerate,eucal,theorem}

\numberwithin{equation}{section}

\theoremheaderfont{\normalfont\scshape}
{\theorembodyfont{\rmfamily}}
{\theorembodyfont{\rmfamily}}
\newtheorem{prop}{Proposition}[section]
\newtheorem{thm}{Theorem}[section]

\DeclareMathOperator{\Aut}{Aut}
\DeclareMathOperator{\Hom}{Hom}
\DeclareMathOperator{\Sing}{Sing}
\DeclareMathOperator{\Sk}{Sk}
\DeclareMathOperator{\Spec}{Spec}

\def\fq{\mathbb{F}_q}
\def\bn{\bigskip\noindent} 
\def\sn{\smallskip\noindent} 

\title{Minimal smooth models of Drinfeld coarse moduli threefolds}
\author{Igor Yu. Potemine}
\date{}

\begin{document}
\maketitle

\begin{abstract}We construct minimal smooth models of coarse moduli 
schemes of rational Drinfeld modules of rank 4 as well as their 
smooth compactifications.
\end{abstract}

\section*{Introduction}

Moduli spaces of Drinfeld modules and shtukas regained a great interest after
the beautiful Lafforgue's proof of the global Langlands conjecture for
$\mathrm{GL}_n$ over global function fields (\emph{cf.}~\cite{La}). The 
smoothness of toroidal Lafforgue compactification of the moduli stack of 
shtukas was crucial in his proof.

\smallskip In this paper we study explicitly the ``toy problem'' of the 
smooth compactification of coarse moduli spaces $M^4(1)$ of rational 
Drinfeld modules of rank $r=4$ (\emph{cf.}~\cite{Po} for the case $r=3$). 
This moduli space is a toric threefold over $A=\fq[T]$. In section 3, 
we construct explicitly its minimal terminal model $M_{\textrm{min}}^4(1)$
which is not smooth. However, it is possible to construct its 
desingularization $M_{\textrm{ess}}^4(1)$ by essential divisors which is 
minimal in the sense of minimal volume of its shed. Finally, in section 5, 
we construct explicitly the minimal terminal and smooth compactifications 
$\overline{M_{\textrm{min}}^4}(1)$ and $\overline{M_{\textrm{ess}}^4}(1)$ 
of our course moduli scheme.

\smallskip The results of this paper are also useful for the theory of
complex multiplication of Drinfeld modules of higher rank 
\cite[ch.~3]{Po:phd}. 

\section{$J$-invariants of Drinfeld modules of rank 4}

Let $A=\fq[T]$ be the ring of polynomials over the finite field $\fq$
and $K=\fq(T)$ its quotient field. Let $L$ be a field equipped with a 
non-trivial morphism $\alpha_L^{}:A\rightarrow L$ and $E$ a Drinfeld module 
of rank 4 over $L$ :
\begin{equation}
T\mapsto T_E^{}=\alpha_L^{}(T)+a_1^{}\tau+a_2^{}\tau^2+a_3^{}\tau^3+
\Delta(E)\tau^4,\quad\Delta(E)\in L^*.
\end{equation}
A coefficient $a_k^{}$, $1\leqslant k\leqslant 4$, is of weight $q^k-1$. The
$\mathbf{j}$-invariant is a triple:
\begin{equation}
\mathbf{j}(E)=(j_1^{}(E),j_2^{}(E),j_3^{}(E))=\biggr(\frac{a_1^{q^3+q^2+q+1}}
{\Delta(E)},\frac{a_2^{q^2+1}}{\Delta(E)},\frac{a_3^{q^3+q^2+q+1}}
{\Delta(E)^{q^2+q+1}}\biggr).
\end{equation}
In addition, we consider the following invariants (necessarily of weight zero)
of isomorphism classes of Drinfeld modules:
\begin{align}
j_{12}^{\delta_1\delta_2}=\frac{a_1^{\delta_1}a_2^{\delta_2}}
{\Delta(E)^{\delta_4}},\quad &\delta_1+\delta_2(q+1)=\delta_4(q^3+q^2+q+1),\\
j_{13}^{\delta_1\delta_3}=\frac{a_1^{\delta_1}a_3^{\delta_3}}
{\Delta(E)^{\delta_4}},\quad &\delta_1+\delta_3(q^2+q+1)=
\delta_4(q^3+q^2+q+1),\\
j_{23}^{\delta_2\delta_3}=\frac{a_2^{\delta_2}a_3^{\delta_3}}
{\Delta(E)^{\delta_4}},\quad &\delta_2(q+1)+\delta_3(q^2+q+1)=
\delta_4(q^3+q^2+q+1),\\
\intertext{and}
j_{123}^{\delta_1\delta_2\delta_3}=\frac{a_1^{\delta_1}a_2^{\delta_2}
a_3^{\delta_3}}{\Delta(E)^{\delta_4}},\quad &\delta_1+\delta_2(q+1)+
\delta_3(q^2+q+1)=\delta_4\frac{q^4-1}{q-1}.
\end{align}
These invariants are called \emph{basic} if 
\begin{equation}
0\leqslant\delta_1,\delta_3\leqslant q^3+q^2+q+1,\textrm{ and } 0\leqslant
\delta_2\leqslant q^2+1.
\end{equation}

\section{Drinfeld coarse moduli threefold}

\medskip Denote $M^4(1)$ the coarse moduli scheme of Drinfeld modules of rank
4. In our previous paper we proved the following result.
\begin{thm}\emph{\cite[th.~3.1]{Po}} The coarse moduli scheme $M^4(1)$ is an 
affine toric $A$-threefold generated by basic $j$-invariants, that is,
\begin{equation}
M^4(1)=\Spec A\big[j_1^{},j_2^{},j_3^{},j_{12}^{\delta_1\delta_2},
j_{13}^{\delta_1\delta_3},j_{23}^{\delta_2\delta_3},
j_{123}^{\delta_1\delta_2\delta_3}\big].
\end{equation}
The $\mathbf{j}$-invariant defines a finite flat covering:
\begin{equation}
M^4(1)\rightarrow \Spec A[j_1^{},j_2^{},j_3^{}]
\end{equation}
of degree $d=(q^3+q^2+q+1)(q^2+1)$.
\end{thm}

This is a singular variety. Its singular locus can be easily described 
(\emph{cf.}~\cite[ex.~4.7]{Po}). Indeed, denote $M^4(1)[2]\subset M^4(1)$ 
is the coarse moduli subscheme of Drinfeld modules of type
\begin{equation}
T\mapsto T_E^{}=\alpha_L^{}(T)+a_2^{}\tau^2+\Delta(E)\tau^4,
\end{equation}
that is, such that $a_1^{}=a_3^{}=0$. 
\begin{prop}\emph{(\cite[th.~3.2]{Po})}
The singular locus of $M^4(1)$ is the affine line
\begin{equation}
\Sing(M^4(1))=M^4(1)[2]=\Spec A[j_2^{}]
\end{equation}
generated by the $j_2^{}$-invariant $j_2^{}(E)=\frac{a_2^{q^2+1}}{\Delta(E)}$.
\end{prop}

It follows from the fact that $\Aut(E)=\mathbb{F}_{q^2}^*/\mathbb{F}_q^*$ is 
non-trivial if $j_1^{}(E)=j_3^{}(E)=0$ and the theorem on 
``purity of branch locus''.

\bigskip So $M^4(1)$ is an affine toric $3$-fold and, therefore, can 
be described by a $3$-dimensional rational polyhedral cone \cite{Da}.  Namely,
we fix a lattice $N^3$ of rank $3$ and let 
$N^{3*}=\Hom_{\mathbb{Z}}(N^3,\mathbb{Z})$ be its dual. There exists a 
natural correspondence between $3$-dimensional rational strictly convex 
polyhedral cones in $N_{\mathbb{R}}^{3*}$ and $3$-dimensional affine toric 
varieties \cite{Da,Fu,Oda}.

$$
\begin{picture}(0,0)%
\includegraphics{dr1.pstex}%
\end{picture}%
\setlength{\unitlength}{3315sp}%
\begingroup\makeatletter\ifx\SetFigFont\undefined%
\gdef\SetFigFont#1#2#3#4#5{%
  \reset@font\fontsize{#1}{#2pt}%
  \fontfamily{#3}\fontseries{#4}\fontshape{#5}%
  \selectfont}%
\fi\endgroup%
\begin{picture}(3952,5338)(2146,-8183)
\put(4891,-8071){\makebox(0,0)[lb]{\smash{\SetFigFont{10}{12.0}{\rmdefault}{\mddefault}{\updefault}$(1,0,0)$}}}
\put(4186,-3001){\makebox(0,0)[lb]{\smash{\SetFigFont{10}{12.0}{\rmdefault}{\mddefault}{\updefault}$(q^2+q+1,0,q^3+q^2+q+1)$}}}
\put(2146,-7051){\makebox(0,0)[lb]{\smash{\SetFigFont{10}{12.0}{\rmdefault}{\mddefault}{\updefault}$(1,q^2+1,0)$}}}
\end{picture}
$$
\begin{figure}[h]
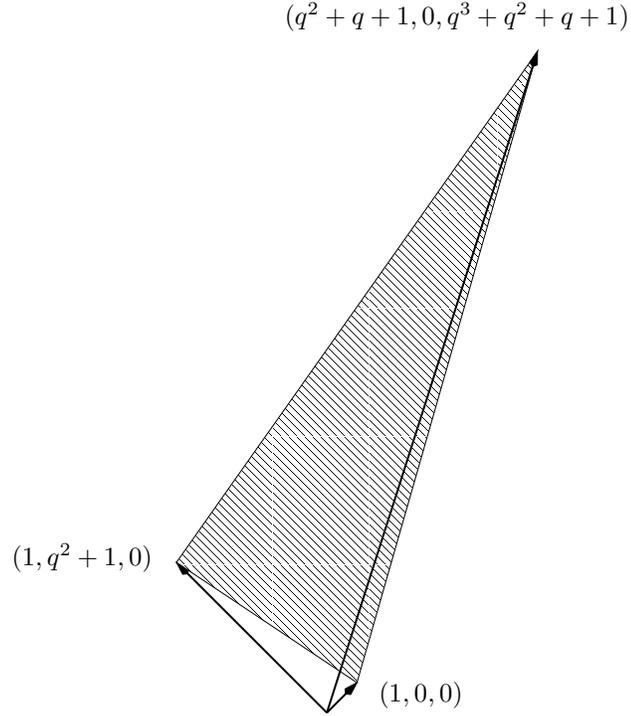

\caption{Dual rational cone $\check\sigma_M^{}$ of $M^4(1)$}{}
\end{figure}

\begin{prop}\emph{\cite[th.~5.1]{Po}}\label{prop:cones} The rational 
simplicial cone $\sigma_M^{}$ of $M^4(1)$ is span\-ned by :
\begin{align}\label{eq:sigmaM}
e_1 &=\left(\frac{q^4-1}{q-1},-q-1,-q^2-q-1\right),\\
e_2 &=(0,1,0),\quad e_3=(0,0,1).
\end{align}
The cone spanned by
\begin{align}
e_1^* &=(1,0,0),\ e_2^*=(1,q^2+1,0),\\
e_3^* &=\left(q^2+q+1,0,\frac{q^4-1}{q-1}\right)
\end{align}
is the dual rational cone $\check\sigma_M^{}$ of $M^4(1)$.
\end{prop}

$$
\begin{picture}(0,0)%
\includegraphics{dr2.pstex}%
\end{picture}%
\setlength{\unitlength}{3315sp}%
\begingroup\makeatletter\ifx\SetFigFont\undefined%
\gdef\SetFigFont#1#2#3#4#5{%
  \reset@font\fontsize{#1}{#2pt}%
  \fontfamily{#3}\fontseries{#4}\fontshape{#5}%
  \selectfont}%
\fi\endgroup%
\begin{picture}(4992,2892)(1328,-3703)
\put(1943,-1006){\makebox(0,0)[lb]{\smash{\SetFigFont{10}{12.0}{\rmdefault}{\mddefault}{\updefault}$(0,0,1)$}}}
\put(1328,-1358){\makebox(0,0)[lb]{\smash{\SetFigFont{10}{12.0}{\rmdefault}{\mddefault}{\updefault}$(0,1,0)$}}}
\put(2820,-3645){\makebox(0,0)[lb]{\smash{\SetFigFont{10}{12.0}{\rmdefault}{\mddefault}{\updefault}$(q^3+q^2+q+1,-q-1,-q^2-q-1)$}}}
\end{picture}
$$
\begin{figure}[h]
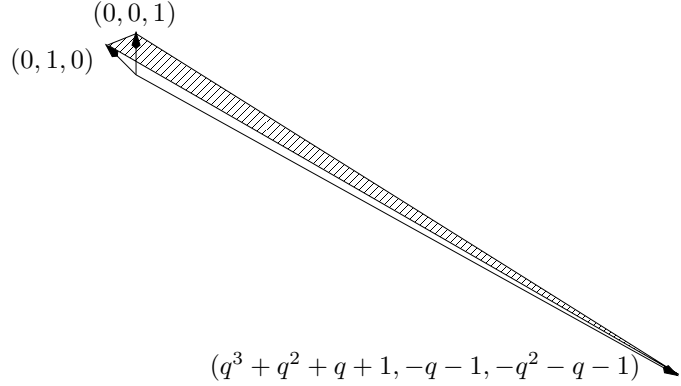

\caption{Rational simplicial cone $\sigma_M^{}$ of $M^4(1)$}{}
\label{fig:cone}
\end{figure}

\section{Minimal terminal model}

In this section we construct the \emph{minimal terminal model} 
$M_{\textrm{min}}^4(1)$ of $M^4(1)$ in the sense of Mori theory applied to the
toric geometry \cite{Re}.

\smallskip We denote $\Sk^1\sigma$ the set of extremal rays of a simplicial 
cone $\sigma$ and $l_\sigma^{}$ a linear form on $N_{\mathbb{Q}}^3$ such that
$l_\sigma^{}\big(\Sk^1\sigma\big)=1$. The convex polytope 
$\sigma\cap l_\sigma^{-1}[0,1]$ is called the \emph{shed} of $\sigma$ and the 
convex polytope $\sigma\cap l_\sigma^{-1}(1)$ in codimension 1 is called the
\emph{roof of the shed} of $\sigma$ (\emph{cf.}~\cite{Re,BGS}). Let now 
$\sigma=\sigma_M^{}$ be the simplicial cone of $M^4(1)$. Then 
$\Sk^1\sigma_M^{}=\{\langle e_1^{}\rangle,\langle e_2^{}\rangle,
\langle e_3^{}\rangle\}$ where $\langle e_i^{}\rangle$, 
$1\leqslant i\leqslant 3$, are rays spanned by vectors $e_i^{}$ defined in
prop.~\ref{prop:cones}. The shed of $\sigma_M^{}$ is the tetrahedron
generated by these vectors and its roof is the facet not containing the origin
(hatched on fig.~\ref{fig:cone}).

\smallskip For any indices $1\leqslant i_1^{}<i_2^{}\leqslant 4$ we will use
the notation $M^4(1)[i_1^{},i_2^{}]$ for the moduli surface of Drinfeld 
modules such that their coefficient $a_i^{}$ ($1\leqslant i\leqslant 3$, 
$i\ne i_1^{},i_2^{}$) is zero. For instance, $M^4(1)[1,3]$ is the moduli
surface of Drinfeld modules such that $a_2^{}=0$.  

\begin{thm}\emph{(\emph{cf.}~\cite[th.~6.1]{Po})}\label{thm:term}
The consecutive star subdivisions centered in the rays 
\begin{equation}
(q^2+1,-1,-q),\quad \textrm{and}\quad (1,0,0) 
\end{equation}
define the unique minimal terminal model $M^4_{\textrm{min}}(1)$ of $M^4(1)$.
\end{thm}

$$
\begin{picture}(0,0)%
\includegraphics{dr3.pstex}%
\end{picture}%
\setlength{\unitlength}{2072sp}%
\begingroup\makeatletter\ifx\SetFigFont\undefined%
\gdef\SetFigFont#1#2#3#4#5{%
  \reset@font\fontsize{#1}{#2pt}%
  \fontfamily{#3}\fontseries{#4}\fontshape{#5}%
  \selectfont}%
\fi\endgroup%
\begin{picture}(4965,6201)(1966,-7177)
\put(6931,-3736){\makebox(0,0)[lb]{\smash{\SetFigFont{9}{10.8}{\rmdefault}{\mddefault}{\updefault}$(0,1,0)$}}}
\put(4171,-1246){\makebox(0,0)[lb]{\smash{\SetFigFont{9}{10.8}{\rmdefault}{\mddefault}{\updefault}$(0,0,1)$}}}
\put(1966,-7096){\makebox(0,0)[lb]{\smash{\SetFigFont{9}{10.8}{\rmdefault}{\mddefault}{\updefault}$(q^3+q^2+q+1,-q-1,-q^2-q-1)$}}}
\put(4591,-3526){\makebox(0,0)[lb]{\smash{\SetFigFont{9}{10.8}{\rmdefault}{\mddefault}{\updefault}$(1,0,0)$}}}
\put(4141,-4786){\makebox(0,0)[lb]{\smash{\SetFigFont{9}{10.8}{\rmdefault}{\mddefault}{\updefault}$(q^2+1,-1,-q)$}}}
\end{picture}
$$
\begin{figure}[h]
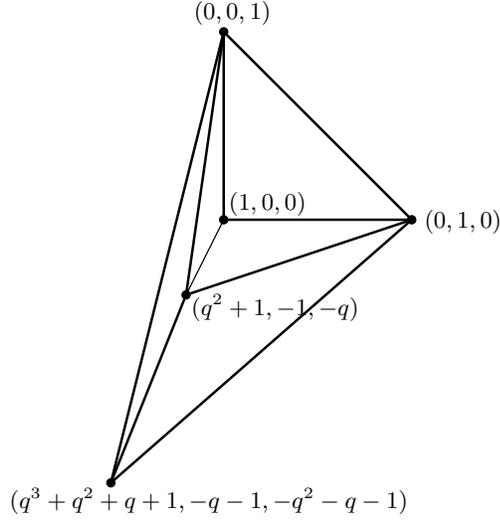

\caption{Rational fan $\Sigma_{\textrm{min}}$ of minimal terminal model
$M^4_{\textrm{min}}(1)$}{}
\label{fig:terminal}
\end{figure}

\sn\emph{Proof.} Denote $\Sigma_{\textrm{min}}$ the fan of 
$M^4_{\textrm{min}}(1)$. The extremal rays of the cones of $\Sigma_{\rm min}$
are called terminal rays of $\sigma_M^{}$. A point of the shed of 
$\sigma_M^{}$ generating a terminal ray will be called a terminal point.

\begin{figure}[h]
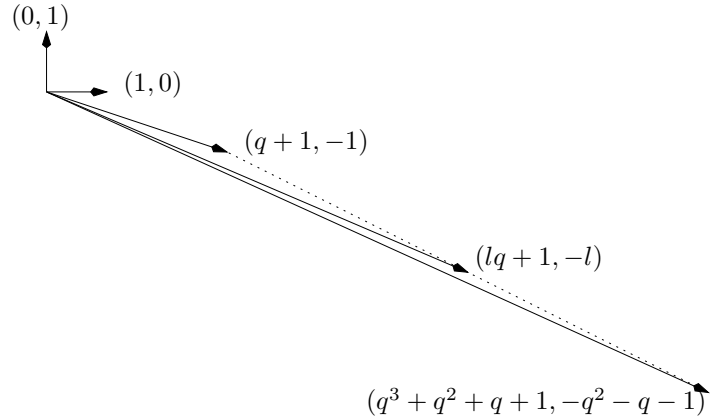

\caption{Minimal desingularization of the surface $M^4(1)[1,3]$}{}
\label{fig:desing}
\end{figure}

Terminal points are necessarily integral points lying inside of the shed of
$\sigma$. Their coordinates may be found by consecutive projections to the
coordinates $(x_1^{},x_3^{})$ and $(-x_3^{},x_2^{})$. The projection on 
$(x_1^{},x_3^{})$ defines the two-dimensional cone
\begin{equation}
\left\langle\left(\frac{q^4-1}{q-1},-q^2-q-1\right),(0,1)\right\rangle.
\end{equation}
which is the rational cone of the surface $M^4(1)[1,3]$. The points 
$(lq+1,-l)$ for $0\leqslant l<q^2+q+1$ are the only integral points
strictly inside of the shed of this cone. They define the minimal 
desingularization of the surface $M^4(1)[1,3]$ 
(\emph{cf.}~fig.\ref{fig:desing}). Note that all these internal terminal 
points belong to the same line $x_1^{}=1-qx_3^{}$.

\begin{figure}[h]
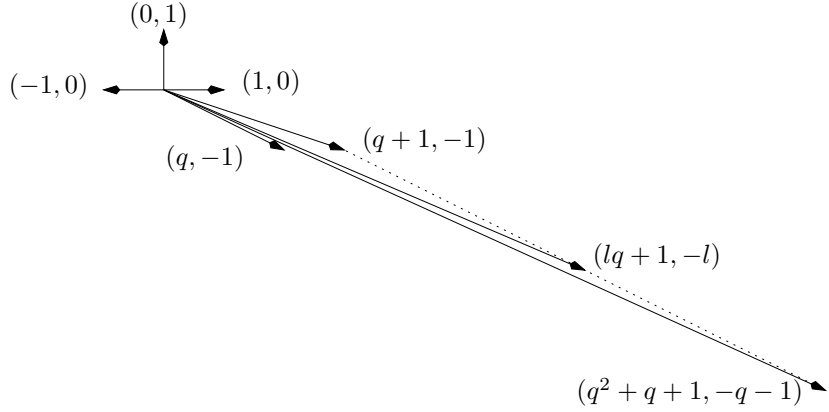

\caption{Minimal smooth compactification of the moduli surface $M^3(1)$}{}
\label{fig:smcomp}
\end{figure}

\smallskip The projection on $(-x_3^{},x_2^{})$ defines the two-dimensional fan
\begin{equation}
\langle(q^2+q+1,-q-1),(0,1),(-1,0)\rangle
\end{equation}
This is the dual fan of the weighted projective space
\begin{equation}
\mathbb{P}_A(q-1,q^2-1,q^3-1)=\mathbb{P}_A(1,q+1,q^2+q+1),
\end{equation}
(over $A$) and this latter can be also considered as the (compactified) 
coarse moduli space $\overline{M^3}(1)$ of Drinfeld modules of rank 
$\leqslant 3$ \cite[sect.~3]{Po}. The points $(lq+1,-l)$ for 
$0\leqslant l<q+1$ and the point $(q,-1)$ are the only integral points 
strictly inside of the shed of this fan. These points give the minimal 
smooth compactification of the moduli surface $M^{3}(1)$
(\emph{cf.}~fig.\ref{fig:smcomp}). Note that all internal terminal 
points (except for $(q,-1)$) belong to the line $x_3^{}=qx_2^{}-1$.

$$
\begin{picture}(0,0)%
\includegraphics{dr4.pstex}%
\end{picture}%
\setlength{\unitlength}{3315sp}%
\begingroup\makeatletter\ifx\SetFigFont\undefined%
\gdef\SetFigFont#1#2#3#4#5{%
  \reset@font\fontsize{#1}{#2pt}%
  \fontfamily{#3}\fontseries{#4}\fontshape{#5}%
  \selectfont}%
\fi\endgroup%
\begin{picture}(5217,3133)(1996,-3839)
\put(1996,-901){\makebox(0,0)[lb]{\smash{\SetFigFont{10}{12.0}{\rmdefault}{\mddefault}{\updefault}$(0,1)$}}}
\put(2836,-1396){\makebox(0,0)[lb]{\smash{\SetFigFont{10}{12.0}{\rmdefault}{\mddefault}{\updefault}$(1,0)$}}}
\put(3736,-1831){\makebox(0,0)[lb]{\smash{\SetFigFont{10}{12.0}{\rmdefault}{\mddefault}{\updefault}$(q+1,-1)$}}}
\put(5461,-2701){\makebox(0,0)[lb]{\smash{\SetFigFont{10}{12.0}{\rmdefault}{\mddefault}{\updefault}$(lq+1,-l)$}}}
\put(4651,-3781){\makebox(0,0)[lb]{\smash{\SetFigFont{10}{12.0}{\rmdefault}{\mddefault}{\updefault}$(q^3+q^2+q+1,-q^2-q-1)$}}}
\end{picture}
$$

\smallskip Combining these results, we get that a point 
$(x_1^{},x_2^{},x_3^{})$ (distinct from the origin and lying strictly inside 
of the shed of $\sigma_M^{}$) can be terminal only if it is $(1,0,0)$, 
$(q^2+1,-1,-q)$ or given by the relations
\begin{equation}\label{eq:terminal}
0\leqslant -x_2^{}\leqslant q,\ x_3^{}=qx_2^{}-1,\textrm{ and }
x_1^{}=1-qx_3^{}.
\end{equation}
On the one hand, an easy straightforward computation shows that points 
\eqref{eq:terminal} lie above the hyperplane passing through 
$e_1^{}$, $e_2^{}$ and $e_3^{}$. So they do not belong to the shed
of $\sigma_M^{}$ and can not be terminal. On the other hand, the rays 
$(q^2+1,-1,-q)$ and $(1,0,0)$ are terminal. The consecutive star 
subdivisions centered in these rays yield 
the shed represented on fig.~\ref{fig:terminal}. This shed has a concave 
roof along internal walls (\emph{cf.}~\cite{Re} for exhaustive terminology).
It follows from Reid's theorem \cite[th.$\:$0.2]{Re} that the corresponding
variety is a minimal terminal model. Any other minimal model with terminal 
singularities should have the same shed. It is easy to check that the roof 
of our shed is strictly concave along internal walls and, consequently, the
constructed minimal model is unique.\hfill$\square$

$$
\begin{picture}(0,0)%
\includegraphics{dr5.pstex}%
\end{picture}%
\setlength{\unitlength}{3315sp}%
\begingroup\makeatletter\ifx\SetFigFont\undefined%
\gdef\SetFigFont#1#2#3#4#5{%
  \reset@font\fontsize{#1}{#2pt}%
  \fontfamily{#3}\fontseries{#4}\fontshape{#5}%
  \selectfont}%
\fi\endgroup%
\begin{picture}(6117,3073)(1096,-3779)
\put(1996,-901){\makebox(0,0)[lb]{\smash{\SetFigFont{10}{12.0}{\rmdefault}{\mddefault}{\updefault}$(0,1)$}}}
\put(2836,-1396){\makebox(0,0)[lb]{\smash{\SetFigFont{10}{12.0}{\rmdefault}{\mddefault}{\updefault}$(1,0)$}}}
\put(3736,-1831){\makebox(0,0)[lb]{\smash{\SetFigFont{10}{12.0}{\rmdefault}{\mddefault}{\updefault}$(q+1,-1)$}}}
\put(5461,-2701){\makebox(0,0)[lb]{\smash{\SetFigFont{10}{12.0}{\rmdefault}{\mddefault}{\updefault}$(lq+1,-l)$}}}
\put(1096,-1426){\makebox(0,0)[lb]{\smash{\SetFigFont{10}{12.0}{\rmdefault}{\mddefault}{\updefault}$(-1,0)$}}}
\put(2266,-1966){\makebox(0,0)[lb]{\smash{\SetFigFont{10}{12.0}{\rmdefault}{\mddefault}{\updefault}$(q,-1)$}}}
\put(5341,-3721){\makebox(0,0)[lb]{\smash{\SetFigFont{10}{12.0}{\rmdefault}{\mddefault}{\updefault}$(q^2+q+1,-q-1)$}}}
\end{picture}
$$

\section{Minimal essential smooth model}

The minimal terminal model $M_{\textrm{min}}^4(1)$ constructed in the 
previous section is not smooth. However, it is possible to construct a unique 
smooth model which is minimal in the sense that the volume (of its shed) is
minimal (\emph{cf.}~\cite{BGS}). Technically speaking, this model is called 
\emph{essential} since it is obtained by a chain of blow-ups of
$M^4(1)$ such that all exceptional divisors are essential, that is, present in
any smooth model \cite{BGS}.

\begin{thm}\label{thm:ess} The star subdivisions of the fan
$\Sigma_{\textrm{min}}$ of $M_{\textrm{min}}^4(1)$ centered in the rays
\begin{equation}\label{eq:ess1}
(kq^2+lq+1,-k,-kq-l),\quad 1\leqslant k\leqslant q,
\end{equation}
followed by the star subdivisions centered in the rays
\begin{equation}\label{eq:ess2}
(kq^2+lq+1,-k,-kq-l),\quad 0\leqslant k\leqslant q,
\end{equation}
where
\begin{equation}
\begin{cases}2\le l\le q &\textrm{if } 1\le k\le q\\
1\le l\le q-1 &\textrm{if }k=0\end{cases},
\end{equation}
define the unique (minimal) essential smooth model $M_{\textrm{ess}}^4(1)$.
\end{thm}

$$
\begin{picture}(0,0)%
\includegraphics{dr6.pstex}%
\end{picture}%
\setlength{\unitlength}{3108sp}%
\begingroup\makeatletter\ifx\SetFigFont\undefined%
\gdef\SetFigFont#1#2#3#4#5{%
  \reset@font\fontsize{#1}{#2pt}%
  \fontfamily{#3}\fontseries{#4}\fontshape{#5}%
  \selectfont}%
\fi\endgroup%
\begin{picture}(10548,9712)(265,-10928)
\put(3601,-1411){\makebox(0,0)[lb]{\smash{\SetFigFont{9}{10.8}{\rmdefault}{\mddefault}{\updefault}$(0,0,1)$}}}
\put(3601,-1861){\makebox(0,0)[lb]{\smash{\SetFigFont{9}{10.8}{\rmdefault}{\mddefault}{\updefault}$(1,0,0)$}}}
\put(6121,-8521){\makebox(0,0)[lb]{\smash{\SetFigFont{9}{10.8}{\rmdefault}{\mddefault}{\updefault}$(0,1,0)$}}}
\put(2701,-10861){\makebox(0,0)[lb]{\smash{\SetFigFont{9}{10.8}{\rmdefault}{\mddefault}{\updefault}$(q^3+q^2+q+1,-q-1,-q^2-q-1)$}}}
\put(4861,-2318){\makebox(0,0)[lb]{\smash{\SetFigFont{9}{10.8}{\rmdefault}{\mddefault}{\updefault}$(q+1,0,-1)$}}}
\put(4981,-2776){\makebox(0,0)[lb]{\smash{\SetFigFont{9}{10.8}{\rmdefault}{\mddefault}{\updefault}$(lq+1,0,-l)$}}}
\put(5086,-3211){\makebox(0,0)[lb]{\smash{\SetFigFont{9}{10.8}{\rmdefault}{\mddefault}{\updefault}$(q^2-q+1,0,1-q)$}}}
\put(519,-8161){\makebox(0,0)[lb]{\smash{\SetFigFont{9}{10.8}{\rmdefault}{\mddefault}{\updefault}$(kq^2+lq+1,-k,-kq-l)$}}}
\put(459,-9060){\makebox(0,0)[lb]{\smash{\SetFigFont{9}{10.8}{\rmdefault}{\mddefault}{\updefault}$(q^3+q+1,-q,-q^2-1)$}}}
\put(265,-9969){\makebox(0,0)[lb]{\smash{\SetFigFont{9}{10.8}{\rmdefault}{\mddefault}{\updefault}$(q^3+lq+1,-q,-q^2-l)$}}}
\put(622,-7269){\makebox(0,0)[lb]{\smash{\SetFigFont{9}{10.8}{\rmdefault}{\mddefault}{\updefault}$(kq^2+q+1,-k,-kq-1)$}}}
\put(892,-6360){\makebox(0,0)[lb]{\smash{\SetFigFont{9}{10.8}{\rmdefault}{\mddefault}{\updefault}$(2q^2+lq+1,-2,-2q-l)$}}}
\put(1156,-5460){\makebox(0,0)[lb]{\smash{\SetFigFont{9}{10.8}{\rmdefault}{\mddefault}{\updefault}$(2q^2+q+1,-2,-2q-1)$}}}
\put(1570,-4568){\makebox(0,0)[lb]{\smash{\SetFigFont{9}{10.8}{\rmdefault}{\mddefault}{\updefault}$(q^2+lq+1,-1,-q-l)$}}}
\put(1855,-3668){\makebox(0,0)[lb]{\smash{\SetFigFont{9}{10.8}{\rmdefault}{\mddefault}{\updefault}$(q^2+q+1,-1,-q-1)$}}}
\put(2357,-3226){\makebox(0,0)[lb]{\smash{\SetFigFont{9}{10.8}{\rmdefault}{\mddefault}{\updefault}$(q^2+1,-1,-q)$}}}
\end{picture}
$$
\begin{figure}[h]
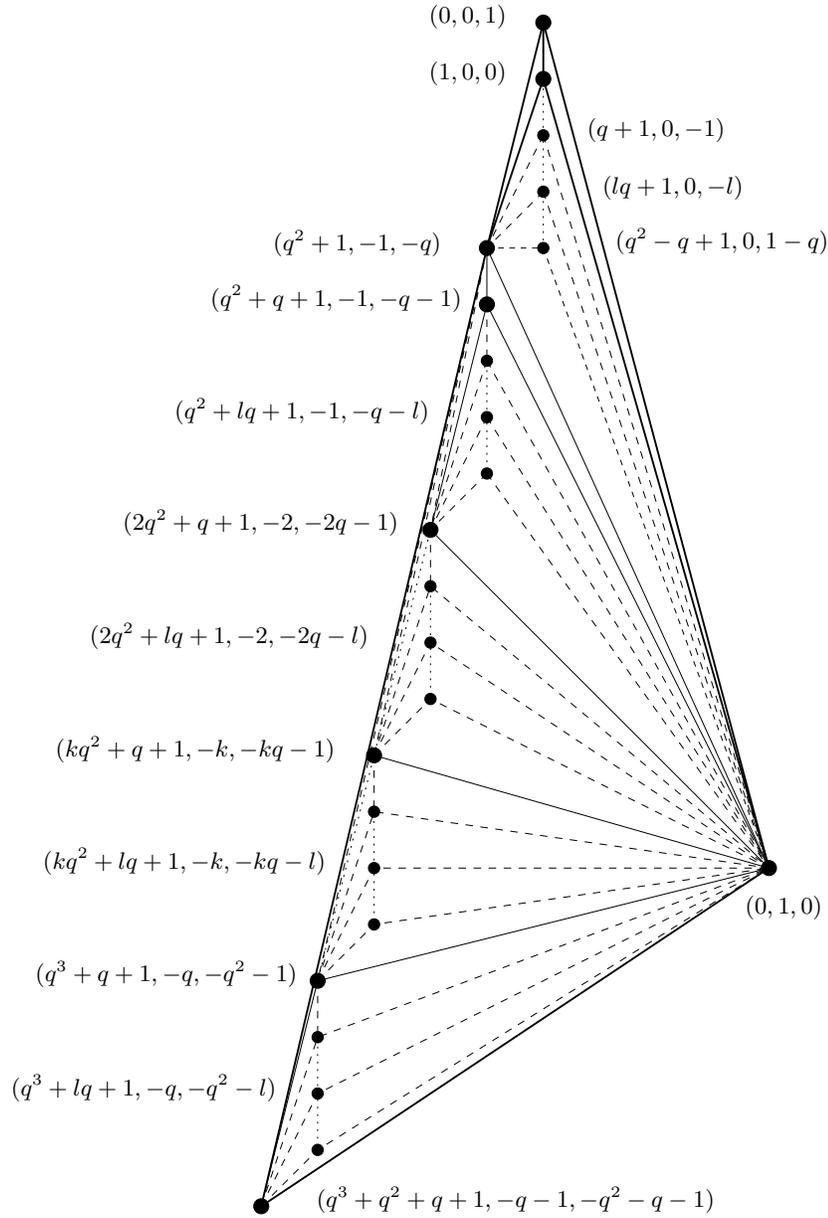

\caption{Minimal essential desingularization of $M^4(1)$}{}
\label{fig:essdesing}
\end{figure}

\sn\emph{Proof.} Any smooth toric $3$-fold corresponds to a \emph{regular}
fan, \emph{i.e.} such that any cone is generated by a basis of $\mathbb{Z}^3$.
Consequently, in order to find a smooth model of $M_{\textrm{min}}^4(1)$, it
is necessary to find a \emph{regular subdivision} of the fan 
$\Sigma_{\textrm{min}}$. Actually, we only need to find a regular 
subdivision of the fan generated by cones
\begin{gather}
\sigma_0'=\left\langle(1,0,0),(0,1,0),(q^2+1,-1,-q)\right\rangle\\
\intertext{and}
\sigma_{q+1}^{}=\left\langle((0,1,0),(q^2+1,-1,-q),
\left(\frac{q^4-1}{q-1},-q-1,-q^2-q-1\right)\right\rangle.
\end{gather}

\smallskip The multiplicity $\mu_{\sigma_{q+1}^{}}^{}$ of the cone
$\sigma_{q+1}^{}$ is equal to $q^2+1$. In view of (\cite{BGS},
Prop.$\;$2.6) there exists a unique point $x_{\sigma_{q+1}^{}}^{}$
of $\sigma_{q+1}^{}$ such that
\begin{equation}
l_{\sigma_{q+1}^{}}(x_{\sigma_{q+1}^{}}^{})=
1+\frac{1}{\mu_{\sigma_{q+1}^{}}^{}}.
\end{equation}
In our case, the linear form $l_{\sigma_{q+1}^{}}^{}$ is given by :
\begin{equation}
l_{\sigma_{q+1}^{}}^{}={{(q^2+2)}\over{q^2+1}}x+y+qz.
\end{equation}
We denote $P_q$ the point $(q^3+q+1,-q,-q^2-1)$ and we obtain that
\begin{equation}
l_{\sigma_{q+1}^{}}^{}(P_q^{})=1+\frac{1}{q^2+1}.
\end{equation}
Using the procedure of $G$-desingularisation \cite{BGS}, we should make
a star subdivision centered in $P_q$.

\smallskip More generally, for any $2\le k\le q+1$, the multiplicity of
the cone 
\begin{equation}
\sigma_k^{}=\langle (0,1,0),(q^2+1,-1,-q),
(kq^2+q+1,-k,-kq-1)\rangle
\end{equation}
is equal to $\mu_{\sigma_k^{}}^{}=(k-1)q+1$. In addition, the point
\begin{equation}
P_{k-1}=\bigr((k-1)q^2+q+1,1-k,(1-k)q-1\bigr)
\end{equation}
is such that
\begin{equation}
l_{\sigma_k^{}}(P_{k-1})=1+\frac{1}{(k-1)q+1}\;.
\end{equation}
In this way we obtain that the consecutive star subdivisions centered
in $P_k$ for $1\leqslant k\leqslant q$ are the first $q$ steps
of the $G$-desingularisation procedure.

\medskip We should now find regular subdivisions of the cones 
\begin{equation}
\sigma_k'=\langle P_k,P_{k+1},(0,1,0)\rangle
\end{equation}
for $1\leqslant k\leqslant q$ and of the cone
\begin{equation}
\sigma_0'=\langle(q^2+1,-1-q),(1,0,0),(0,1,0)\rangle.
\end{equation}
These are cones of multiplicity $q$. It is easy to notice that the points 
\begin{equation}
P_{kl}=(kq^2+lq+1,-k,-kq-l),
\end{equation} 
for any $2\leqslant l\leqslant q$ if $1\leqslant k\leqslant q$ and for any 
$1\leqslant l\leqslant q-1$ if $k=0$, are aligned for a fixed $k$ and such
that 
\begin{equation}
l_{\sigma_k'}(P_{kl})\in [1,2[.
\end{equation}
Consequently, we have obtained the minimal desingularisation of $M^4(1)$ by
essential divisors.\hfill$\square$

\section{Minimal smooth compactificaion}

Denote $\overline{M^4}(1)$ the moduli space of Drinfeld modules of rank
$\leqslant 4$. On the one hand, it is easy to prove that this space is a 
canonical compactification of $M^4(1)$. On the other hand, it is easy to
describe it in terms of toric geometry.

\begin{prop}\emph{(\cite[1.6]{Ka}, \cite[prop.~3.3]{Po})}
The weighted projective space
\begin{align}
\overline{M^4}(1) &=\mathbb{P}_A(q-1,q^2-1,q^3-1,q^4-1)\\
 &=\mathbb{P}_A(1,q+1,q^2+q+1,q^3+q^2+q+1),
\end{align}
is the coarse moduli scheme of rational Drinfeld modules of rank 
$\leqslant 4$. The fan spanned by
\begin{align}
e_1 &=\left(\frac{q^4-1}{q-1},-q-1,-q^2-q-1\right),\\
e_2 &=(0,1,0),\ e_3=(0,0,1)\ \mathrm{and}\ 
e_4^{}=(-1,0,0),
\end{align}
is the rational simplicial fan $\overline{\sigma}_M^{}$ of 
$\overline{M^4}(1)$.
\end{prop}
Indeed, it is easy to see (\emph{cf.}~\cite[prop.~3.3]{Po}) that
the affine subvariety of
\begin{equation}
\overline{M^4}(1)=\mathbb{P}_A(q-1,q^2-1,q^3-1,q^4-1)
\end{equation}
corresponding to the non-zero $k$th coordinate, $1\leqslant k\leqslant 3$, 
is the coarse moduli scheme $\overline{M^4}(1)[j_k^{}\ne0]$ of Drinfeld 
modules of rank $\leqslant r$ with non-zero $j_k^{}$-invariant. Moreover, 
$\overline{M^4}(1)$ is the gluing of $M^4(1)$ with these three affine 
$A$-varieties. Its fan $\overline{\sigma}_M^{}$ is obtained by adding the
ray $(-1,0,0)$ to the rational cone $\sigma_M^{}$ 
(\emph{cf.}~prop.~\ref{prop:cones}).

$$\begin{picture}(0,0)%
\includegraphics{dr7.pstex}%
\end{picture}%
\setlength{\unitlength}{3315sp}%
\begingroup\makeatletter\ifx\SetFigFont\undefined%
\gdef\SetFigFont#1#2#3#4#5{%
  \reset@font\fontsize{#1}{#2pt}%
  \fontfamily{#3}\fontseries{#4}\fontshape{#5}%
  \selectfont}%
\fi\endgroup%
\begin{picture}(5112,2862)(1208,-3703)
\put(2820,-3645){\makebox(0,0)[lb]{\smash{\SetFigFont{10}{12.0}{\rmdefault}{\mddefault}{\updefault}$(q^3+q^2+q+1,-q-1,-q^2-q-1)$}}}
\put(2288,-1036){\makebox(0,0)[lb]{\smash{\SetFigFont{10}{12.0}{\rmdefault}{\mddefault}{\updefault}$(0,0,1)$}}}
\put(1208,-1613){\makebox(0,0)[lb]{\smash{\SetFigFont{10}{12.0}{\rmdefault}{\mddefault}{\updefault}$(-1,0,0)$}}}
\put(1343,-1126){\makebox(0,0)[lb]{\smash{\SetFigFont{10}{12.0}{\rmdefault}{\mddefault}{\updefault}$(0,1,0)$}}}
\end{picture}
$$ 
\begin{figure}[h]
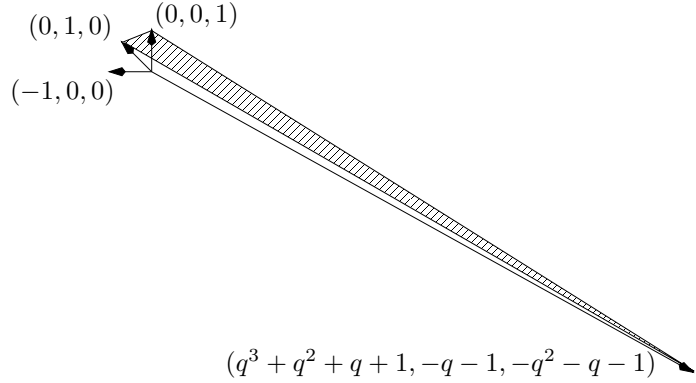

\caption{Rational simplicial fan $\overline{\sigma}_M^{}$ of 
$\overline{M^4}(1)$}{}
\end{figure}

\sn Now we can extend our theorems \ref{thm:term} and \ref{thm:ess}.

\begin{thm} i) The consecutive star subdivisions of $\overline{\sigma}_M^{}$
centered in the rays 
\begin{equation*}
(q^2+1,-1,-q),\quad \textrm{and}\quad (1,0,0) 
\end{equation*}
followed by the star subdivision centered in $(q,0,-1)$ define the unique 
minimal terminal compactification $\overline{M_{\textrm{min}}^4}(1)$ of 
$M^4(1)$.

\sn ii) The consecutive star subdivisions of $\overline{\sigma}_M^{}$ 
centered in the rays \eqref{eq:ess1}, \eqref{eq:ess2} and
\begin{equation}
(kq^2+q,-k,-kq-1)
\end{equation}
for $0\leqslant k\leqslant q$ define the unique minimal 
essential smooth compactification of $M^4(1)$.
\end{thm}

\sn\emph{Proof.} We have seen above that $\overline{M^4}(1)$ is the gluing of
$M^4(1)$ with $\overline{M^4}(1)[j_k^{}\ne0]$ for $1\leqslant k\leqslant 3$.
The affine subvariety $\overline{M^4}(1)[j_1^{}\ne0]$ is isomorphic to 
$\mathbb{A}_A^3$ and, consequently, non-singular. In addition, it is easy to
see that $\overline{M^4}(1)[j_2^{}\ne0]$ corresponding to the cone
\begin{equation}
\langle(-1,0,0),(0,0,1),(q^3+q^2+q+1,-q-1,-q^2-q-1)\rangle
\end{equation}
is desingularized by the blow-up corresponding to the star subdivision 
centered in the ray $(q^2+1,-1,-q)$ belonging to the plane
\begin{equation}
\langle O,(0,0,1),(q^3+q^2+q+1,-q-1,-q^2-q-1)\rangle.
\end{equation}
$$
\begin{picture}(0,0)%
\includegraphics{dr8.pstex}%
\end{picture}%
\setlength{\unitlength}{2901sp}%
\begingroup\makeatletter\ifx\SetFigFont\undefined%
\gdef\SetFigFont#1#2#3#4#5{%
  \reset@font\fontsize{#1}{#2pt}%
  \fontfamily{#3}\fontseries{#4}\fontshape{#5}%
  \selectfont}%
\fi\endgroup%
\begin{picture}(6383,4468)(1531,-5519)
\put(4141,-3571){\makebox(0,0)[lb]{\smash{\SetFigFont{9}{10.8}{\rmdefault}{\mddefault}{\updefault}$O$}}}
\put(4171,-1246){\makebox(0,0)[lb]{\smash{\SetFigFont{9}{10.8}{\rmdefault}{\mddefault}{\updefault}$(0,0,1)$}}}
\put(5626,-3976){\makebox(0,0)[lb]{\smash{\SetFigFont{9}{10.8}{\rmdefault}{\mddefault}{\updefault}$(q^2+1,-1,-q)$}}}
\put(1531,-3526){\makebox(0,0)[lb]{\smash{\SetFigFont{9}{10.8}{\rmdefault}{\mddefault}{\updefault}$(-1,0,0)$}}}
\put(4726,-5461){\makebox(0,0)[lb]{\smash{\SetFigFont{9}{10.8}{\rmdefault}{\mddefault}{\updefault}$(q^3+q^2+q+1,-q-1,-q^2-q-1)$}}}
\end{picture}
$$
\begin{figure}[h]
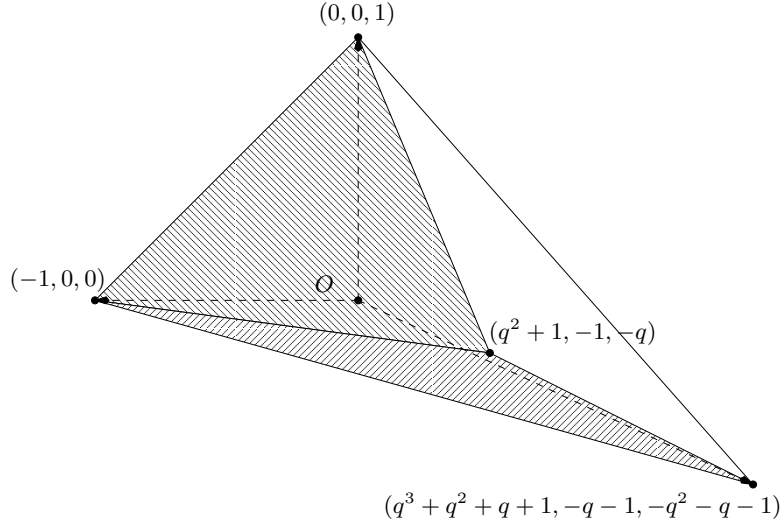

\caption{Desingularization of $\overline{M^4}(1)[j_2^{}\ne0]$}{}
\end{figure}

\smallskip Thus, we only need to find a minimal terminal and essential 
smooth models of the affine $A$-variety $\overline{M^4}(1)[j_3^{}\ne0]$ 
corresponding to the cone
\begin{equation}\label{j3not0}
\langle(-1,0,0),(0,1,0),(q^3+q^2+q+1,-q-1,-q^2-q-1)\rangle
\end{equation}
as well as its desingularization by essential divisors. The proof of i) is 
similar to the proof of theorem~\ref{thm:term} (\emph{cf.}~also 
\cite[th~6.2]{Po}). Indeed, consider projections of the fan 
$\overline{\sigma}_M^{}$ to the coordinates $(x_1^{},x_3^{})$ and 
$(-x_3^{},x_2^{})$. These projections define two-dimensional fans 
\begin{equation}
\left\langle(-1,0),(0,1),\left(\frac{q^{k+1}-1}{q-1},\frac{1-q^k}{q-1}
\right)\right\rangle
\end{equation}
for $k=2,3$. Searching for integral points inside the sheds of these fans, 
one finds two points of theorem~\ref{thm:term} as well as the point $(q,0,-1)$
lying in the cone
\begin{equation}
\left\langle(-1,0),(0,1),\left(q^3+q^2+q+1,-q^2-q-1\right)\right\rangle.
\end{equation}
Star subdivisions centered in corresponding terminal rays define 
$\overline{M_{\textrm{min}}^4}(1)$.

\smallskip In order to prove ii), we should now desingularize  
$\overline{M^4}(1)[j_3^{}\ne0]$ corresponding to the cone
\eqref{j3not0}. The cones
\begin{equation}
\langle(0,1,0),(q,0,-1),(q^3+q^2+q+1,-q-1,-q^2-q-1)\rangle
\end{equation}
and
\begin{equation}
\langle(-1,0,0),(0,1,0),(q,0,-1)\rangle
\end{equation}
are already regular. It suffices to find a regular subdivision of the cone 
\begin{equation}
\widetilde{\sigma}_0^{}=\langle(-1,0,0),(q,0,-1),
(q^3+q^2+q+1,-q-1,-q^2-q-1)\rangle.
\end{equation}
of multiplicity $q+1$. 

$$
\begin{picture}(0,0)%
\includegraphics{dr9.pstex}%
\end{picture}%
\setlength{\unitlength}{3315sp}%
\begingroup\makeatletter\ifx\SetFigFont\undefined%
\gdef\SetFigFont#1#2#3#4#5{%
  \reset@font\fontsize{#1}{#2pt}%
  \fontfamily{#3}\fontseries{#4}\fontshape{#5}%
  \selectfont}%
\fi\endgroup%
\begin{picture}(2925,4618)(2161,-5204)
\put(5086,-1411){\makebox(0,0)[lb]{\smash{\SetFigFont{10}{12.0}{\rmdefault}{\mddefault}{\updefault}$(0,1,0)$}}}
\put(4006,-781){\makebox(0,0)[lb]{\smash{\SetFigFont{10}{12.0}{\rmdefault}{\mddefault}{\updefault}$(-1,0,0)$}}}
\put(3466,-1861){\makebox(0,0)[lb]{\smash{\SetFigFont{10}{12.0}{\rmdefault}{\mddefault}{\updefault}$(q,0,-1)$}}}
\put(2161,-3616){\makebox(0,0)[lb]{\smash{\SetFigFont{10}{12.0}{\rmdefault}{\mddefault}{\updefault}$(q^3+q,-q,-q^2-1)$}}}
\put(2566,-5146){\makebox(0,0)[lb]{\smash{\SetFigFont{10}{12.0}{\rmdefault}{\mddefault}{\updefault}$(q^3+q^2+q+1,-q-1,-q^2-q-1)$}}}
\put(2206,-2761){\makebox(0,0)[lb]{\smash{\SetFigFont{10}{12.0}{\rmdefault}{\mddefault}{\updefault}$(kq^2+q,-k,-kq-1)$}}}
\end{picture}
$$
\begin{figure}[h]
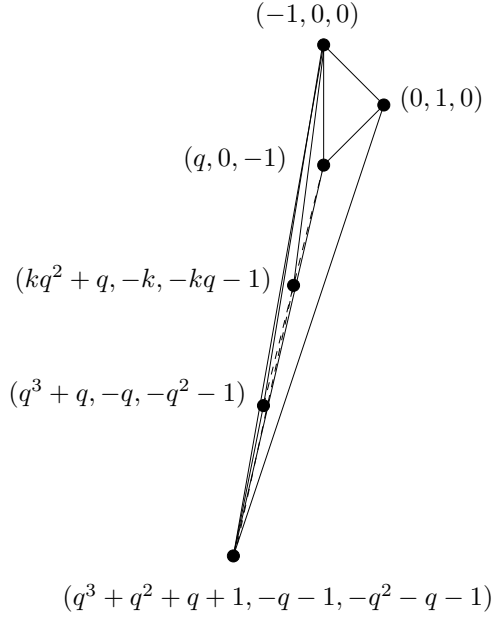

\caption{Desingularization of $\overline{M^4}(1)[j_3^{}\ne0]$}{}
\end{figure}

\smallskip The linear form $l_{\widetilde{\sigma}_0^{}}^{}$ is given 
by :
\begin{equation}
l_{\widetilde{\sigma}_0^{}}^{}=-x+\frac{q^2+q+1}{q+1}y-(q+1)z.
\end{equation}
Denote $Q_1$ the point $(q^2+q,-1,-q-1)$ then we obtain that
\begin{equation}
l_{\widetilde{\sigma}_0^{}}^{}(Q_1^{})=1+\frac{1}{q+1}.
\end{equation}
Using the procedure of $G$-desingularisation \cite{BGS}, we should make
a star subdivision centered in $Q_1$. More generally, for any $0\le k\le q-1$,
the multiplicity of the cone 
\begin{equation}
\widetilde{\sigma}_k^{}=\left\langle (-1,0,0),(kq^2+q,-k,-kq-1),
\left(\frac{q^4-1}{q-1},-q-1,\frac{1-q^3}{q-1}\right)\right\rangle
\end{equation}
is equal to $\mu_{\widetilde{\sigma}_k^{}}^{}=q+1-k$. In addition, the point
\begin{equation}
Q_{k+1}=\bigr((k+1)q^2+q,1-(k+1),-(k+1)q-1\bigr)
\end{equation}
is such that
\begin{equation}
l_{\sigma_k^{}}(Q_{k+1})=1+\frac{1}{q+1-k}\;.
\end{equation}
In this way we obtain that the consecutive star subdivisions centered
in $Q_k$ for $1\leqslant k\leqslant q$ define the $G$-desingularisation 
of $\widetilde{\sigma}_0{}$.

\smallskip Thus, we have found the minimal desingularization of 
$\overline{M^4}(1)[j_3^{}\ne0]$ and it finishes the proof.\hfill$\square$

\bn\textbf{Acknowledgments.} I am very thankful to Marc
Reversat for his warm and constant support.

\nocite{*}
\bibliographystyle{amsplain}
\bibliography{dr3folds}

\bigskip
\begin{flushright}
Igor Potemine\\
Laboratoire Emile Picard\\
Universit\'e Paul Sabatier\\
118, route de Narbonne\\
31062 Toulouse C\'edex 4\\
France
\end{flushright}

\begin{flushright}
e-mail : potemine@picard.ups-tlse.fr
\end{flushright}

\end{document}